\documentclass{article}

\usepackage[latin1]{inputenc}
\usepackage{graphicx}
\usepackage{amsmath}
\usepackage{amssymb}
\usepackage{amsthm}
\usepackage{amsfonts}
\usepackage{amstext}
\usepackage{amsopn}
\usepackage{amsxtra}
\usepackage{mathrsfs}
\usepackage{fourier}
\usepackage{dsfont}
\usepackage{esint}
\usepackage{enumitem}
\usepackage{dsfont}
\usepackage[colorlinks=true]{hyperref}
\hypersetup{urlcolor=blue, citecolor=red, linkcolor=blue}
\usepackage{bookmark}
\usepackage[usenames,dvipsnames]{color}
\usepackage{setspace}

\newtheorem{thm}{Theorem}
\newtheorem{lemma}[thm]{Lemma}
\newtheorem{prop}[thm]{Proposition}
\newtheorem{cor}[thm]{Corollary}

\newcounter{AMSsection}
\newcommand{\AMSsection}[1]{
	\setcounter{AMSsubsection}{0}\refstepcounter{AMSsection}
	\begin{center}
		{\sc \theAMSsection.~#1}\end{center}}
\newcounter{AMSsubsection}
\newcommand{\AMSsubsection}[1]{
	\stepcounter{AMSsubsection}
	\medskip\noindent\theAMSsection.\theAMSsubsection.~{\bf #1}.~}
\newcommand{\labelSec}[1]{\hypertarget{#1}{\kern 0.1pt}\label{#1}}
\newcommand{\refSec}[1]{\hyperlink{#1}{\ref{#1}}}

\newcommand{\Email}[1]{{\footnotesize{{\sl E-mail:\/} \href{mailto:#1}{\rm\textsf{#1}}}}}
\newcommand{\be}[1]{\begin{equation}\label{#1}}
\newcommand{\ee}{\end{equation}}
\renewcommand{\(}{\left(}
\renewcommand{\)}{\right)}
\usepackage{color}

\newcommand{\bq}{\begin{equation}}

\newcommand{\DD}{\mathsf D}

\newcommand{\dt}{d\kern-0.5pt t}

\newcommand{\idB}[2]{\int_{\partial B_{#1}}#2\,d\kern0.5pt\varsigma}

\newcommand{\ird}[1]{\int_{\R^d}{#1}\,dx}
\newcommand{\irdmu}[1]{\int_{\R^d}#1\,d\mu}
\newcommand{\iS}[1]{\int_{\S^d}{#1}\,d\mu}

\newcommand{\Lap}{\Delta}

\newcommand{\nrm}[2]{\|{#1}\|_{\mathrm L^{#2}(\S^d)}}
\newcommand{\nrmC}[2]{\|{#1}\|_{\mathrm L^{#2}(\mathcal C)}}
\newcommand{\nrmcnd}[2]{\|{#1}\|_{\mathrm L^{#2}(\mathcal C)}}

\newcommand{\nrmR}[2]{\|{#1}\|_{\mathrm L^{#2}(\R)}}
\newcommand{\nrmRd}[2]{\|{#1}\|_{\mathrm L^{#2}(\R^d)}}

\newcommand{\R}{\mathbb R}

\newcommand{\Sp}{\mathbb S}

\newcommand{\courbe}{\mathcal B}
\newcommand{\iM}[1]{\int_{\mathcal M}{#1}\,dv_g}

\renewcommand{\L}{\mathcal L}

\renewcommand{\S}{{\mathbb S}}
\renewcommand{\S}{\mathbb{S}}

\renewcommand\phi{\varphi}

\begin{document}

\title{Symmetry and symmetry breaking:\\ rigidity and flows in elliptic PDEs}

\author{Jean Dolbeault, Maria J. Esteban, and Michael Loss}

\date{}

\maketitle
\thispagestyle{empty}

\begin{abstract}
The issue of symmetry and symmetry breaking is fundamental in all areas of science. Symmetry is often assimilated to order and beauty while symmetry breaking is the source of many interesting phenomena such as phase transitions, instabilities, segregation, self-organization, \emph{etc.} In this contribution we review a series of sharp results of symmetry of nonnegative solutions of nonlinear elliptic differential equation associated with minimization problems on Euclidean spaces or manifolds. Nonnegative solutions of those equations are unique, a property that can also be interpreted as a rigidity result. The method relies on linear and nonlinear flows which reveal deep and robust properties of a large class of variational problems. Local results on linear instability leading to symmetry breaking and the bifurcation of non-symmetric branches of solutions are reinterpreted in a larger, global, variational picture in which our flows characterize directions of descent.
\end{abstract}

\renewcommand{\thefootnote}{\*}\footnote{\hspace*{-14pt}{\small
\emph{Acknowledgements:} Partially supported by the projects \emph{Kibord} and \emph{EFI} (J.D.) of the French National Research Agency (ANR), and by the NSF grant DMS-1301555 (M.L.).
}}

\AMSsection{Introduction}\labelSec{Sec:Intro}

\emph{Symmetries} are fundamental properties of the laws of Physics. They impose constraints on modeling phenomena and, at a more basic level, they serve as criteria of classification. Inspired by his work in crystallography, Pierre Curie made an early attempt (in 1894) to investigate the consequences of symmetries. Since then, symmetry has been an important preoccupation for many scientists.

More intriguing than symmetry is the phenomenon of \emph{symmetry breaking}, which asserts that the state of a system may have less symmetries than the underlying physical laws. Among various considerations on the causes of the symmetries and what these symmetries mean in physics, P.~Curie wrote in~\cite{curie:jpa-00239814} that 
\begin{center}
\emph{C'est la dissym\'etrie qui cr\'ee le ph\'enom\`ene.}
\end{center}
{}In mathematical terms, ``dissym\'etrie'' shifts the attention to solutions which may have less symmetries than the problem they solve. \emph{Symmetry breaking}, especially spontaneous symmetry breaking, has been an incredibly fruitful concept over the last century. It appears in mechanics (buckling instabilities), in particle physics, in the description of phase transitions or complex dynamics, \emph{etc.} One of the basic mechanisms is the bifurcation phenomenon in nonlinear systems, which has to do with the stability analysis of symmetric states. 

\smallskip Symmetry has attracted the attention of mathematicians for diverse reasons which range from assertions like ``symmetry is beautiful'' to practical motivations: symmetry simplifies the search of solutions and makes their computation more tractable from a numerical point of view by reducing the number of degrees of freedom. 

\emph{Entropy methods} have a long history in various fields of Science and in particular of Mathematics. The notion of entropy that we shall consider here is inspired by results in the theory of nonlinear PDEs and especially nonlinear diffusion equations. It borrows tools from Kinetic Theory and from Information Theory. Other major sources of inspiration are the \emph{carr\'e du champ} method used in the study of Semi-groups and Markov processes as well as the \emph{rigidity} (uniqueness) techniques in the Theory of Nonlinear Elliptic Equations. In addition to the application to symmetry issues, one of our contributions was to rephrase these two approaches in a common framework of parabolic equations and to emphasize the role of the nonlinear diffusions in the search for optimal ranges and optimal constants in related interpolation inequalities. 

It is definitely out of reach to give even a partial account of all mathematical issues of symmetry and symmetry breaking in this paper, so we shall focus on PDEs with two main examples: the first one is the equation
\[
-\,\mbox{div}\,\big(|x|^{-\beta}\,\nabla w\big)=|x|^{-\gamma}\,\big(w^{2p-1}-w^p\big)\quad\mbox{in}\quad\R^d\setminus\{0\}\,,
\]
which has an interesting feature: there is a competition between nonlinearities and weights.
The solutions can be interpreted as critical points of an energy functional.
Without weights, solutions are radially symmetric (up to translations). With weights and
in some regime of the parameters $\beta$, $\gamma$ and $p$, non-radial solutions are energetically more favorable. Since we are interested in energy minimizers, as a particular sub-problem, understanding who wins in the competition is a central question.

Alternatively, we shall consider the equation
\[
-\Delta\,\varphi+\Lambda\,\varphi=\varphi^{p-1}\quad\mbox{on}\quad\mathcal M\,,
\]
where $\mathcal M$ is a sphere, a compact manifold or a cylinder. In that case, the geometric properties of the manifold replace the weight and compete with the scale induced by the parameter $\Lambda$. If there is enough space, in a precise sense that can be measured, then solutions with less symmetry may have a lower energy.

These two equations, although very simple because the nonlinearities (and also the weights in the case of the first equation) obey power laws, are not purely academic. For one, the solutions (and the associated functional inequalities) are of direct interest for instance in some models of fluid mechanics. More important is the fact that power laws appear in many problems when scalings or blow-up methods are used to extract an asymptotic behavior. Hence, we expect that our model equations lie at the core of many nonlinear or weighted problems. Finally, models involving power laws have the advantage that they can be treated by using \emph{nonlinear flows} and \emph{entropy methods}. Indeed we are able to give sharp results of rigidity for the equation, and symmetry results for the optimal functions associated with related interpolation inequalities. 

Because of the confluence of various branches of analysis such as non-linear diffusion and the calculus of variations, and the fundamental nature of the above equations, we believe that it is worth studying them in great detail, with sharp stability results and sharp constants in the functional inequalities. Note that this amounts to establishing the exact range of the parameters for which extremal functions are symmetric. Variational issues of the symmetry and symmetry breaking will be detailed below.

\smallskip Let us fix some notations and conventions. Throughout this paper, we shall use the notation $2^*:=\frac{2\,d}{d-2}$ if $d\ge3$, and $2^*:=\infty$ if $d=1$ or $2$. We shall say that a function is an \emph{extremal function} for an optimal functional inequality if equality holds in the inequality. To simplify notations, parameters will be omitted whenever they are not essential for the understanding of the strategy of proof. This paper is a review of various results which were published in several papers (references will appear in the text) and are collected together for the first time. The reader is invited to pay attention that some notations have been redefined compared to the original papers.

\AMSsection{Interpolation inequalities and flows on compact manifolds}\labelSec{Sec:sphere}

\AMSsubsection{Interpolation inequalities on \texorpdfstring{$\S^d$}{Sd}}\labelSec{Sec:SPhereInterp}
Let us consider the inequality
\be{Ineq:GNS}
\nrm{\nabla u}2^2+\frac d{p-2}\,\nrm u2^2\ge\frac d{p-2}\,\nrm up^2\quad\forall\,u\in\mathrm H^1(\S^d,d\mu)
\ee
where $d\mu$ is the uniform probability measure induced by the Lebesgue measure on $\S^d\subset\R^{d+1}$. Here the exponent $p$ is such that $1\le p<2$ or $2<p<2^*$, or $p=2^*$ if $d\ge3$. 
The case $p=2^*$ corresponds to the usual Sobolev inequality on $\S^d$ or, using the stereographic projection, to the Sobolev inequality in $\R^d$. In the limit case as $p\to2$, we recover the logarithmic Sobolev inequality
\be{Ineq:logSob}
\nrm{\nabla u}2^2\ge\frac d2\,\iS{|u|^2\,\log\(\frac{|u|^2}{\nrm u2^2}\)}\quad\forall\,u\in\mathrm H^1(\S^d,d\mu)\setminus\{0\}\,.
\ee
In~\eqref{Ineq:GNS} and~\eqref{Ineq:logSob}, equality is achieved by any constant non-zero function. The value of the optimal constants, $d/(p-2)$ and $d/2$ is obtained by linearization: if $\varphi$ is an eigenfunction associated with the first positive eigenvalue of the Laplace-Beltrami operator on $\S^d$, the infimum of
\[
\frac{(p-2)\,\nrm{\nabla u}2^2}{\nrm up^2-\nrm u2^2}\quad\mbox{and}\quad\frac{2\,\nrm{\nabla u}2^2}{\iS{|u|^2\,\log\Big(\frac{|u|^2}{\nrm u2^2}\Big)}}\,,
\]
respectively for $p\neq2$ and for $p=2$, is achieved by $u=1+\varepsilon\,\varphi$ in the limit as $\varepsilon\to0$. 

Inequality~\eqref{Ineq:GNS} has been established in~\cite{BV-V} by rigidity methods, in~\cite{MR1230930} by techniques of harmonic analysis, and using the \emph{carr\'e du champ} method in~\cite{MR1231419,MR1412446,MR2381156}, for any $p>2$. The case $p=2$ was studied in~\cite{MR674060}. In~\cite{MR772092,Bakry-Emery85,MR808640}, D.~Bakry and M.~Emery proved the inequalities under the restriction
\[
2<p\le2^\#:=\frac{2\,d^2+1}{(d-1)^2}\,.
\]
Their method relies on a linear heat flow method which is presented below, as well as a nonlinear flow which allow us to get rid of this restriction.

\AMSsubsection{Flows and \emph{carr\'e du champ} methods on \texorpdfstring{$\S^d$}{Sd}}\labelSec{Sec:Flows}
We start by the linear heat flow method of~\cite{MR808640}. For any function $\rho>0$ we define a \emph{generalized entropy functional} $\mathcal E_p$ and a \emph{generalized Fisher information functional} $\mathcal I_p$ by
\[
\mathcal E_p[\rho]:=\frac1{p-2}\,\left[\(\iS\rho\)^\frac2p-\iS{\rho^\frac2p}\right]\;\mbox{and}\;\mathcal E_2[\rho]:=\frac12\iS{\rho\,\log\(\frac\rho{\nrm\rho1}\)}
\]
if $p\neq2$ or $p=2$, respectively, and
\[
\mathcal I_p[\rho]:=\iS{|\nabla\rho^\frac1p|^2}\,.
\]
With this notation,~\eqref{Ineq:GNS} and~\eqref{Ineq:logSob} amount to $\mathcal I_p[\rho]\ge d\,\mathcal E_p[\rho]$ as can be checked using $\rho=|u|^p$. Let us consider the heat flow
\be{Eqn:Heat}
\frac{\partial\rho}{\partial t}=\Delta\rho
\ee
where $\Delta$ denotes the Laplace-Beltrami operator on $\S^d$, and compute
\[
\frac d{dt}\mathcal E_p[\rho]=-\,2\,\mathcal I_p[\rho]\quad\mbox{and}\quad\frac d{dt}\mathcal I_p[\rho]\le-\,2\,d\,\mathcal I_p[\rho]
\]
where the differential inequality holds if $p\le2^\#$. Under this condition, we obtain that
\[
\frac d{dt}\Big(\mathcal I_p[\rho]-\,d\,\mathcal E_p[\rho]\Big)\le0\,.
\]
On the other hand, $\rho(t,\cdot)$ converges as $t\to\infty$ to a constant, namely $\iS\rho$ since $d\mu$ is a probability measure and $\iS\rho$ is conserved by~\eqref{Eqn:Heat}. As a consequence, $\lim_{t\to\infty}\(\mathcal I_p[\rho]-d\,\mathcal E_p[\rho]\)=0$, which proves that $\mathcal I_p[\rho(t,\cdot)]$ $-\,d\,\mathcal E_p[\rho(t,\cdot)]$ is nonnegative for any $t\ge0$ and completes the proof. See~\cite{MR808640} for details. One may wonder whether the monotonicity property is also true for some $p>2^\#$. The following result contains a negative answer to this question.
\begin{prop}\label{Prop:Counter-Example}\cite{DEL-tlse} For any $p\in(2^\#,2^*)$ or $p=2^*$ if $d\ge3$, there exists a function $\rho_0$ such that, if $\rho$ is a solution of~\eqref{Eqn:Heat} with initial datum $\rho_0$, then
\[
\frac d{dt}\Big(\mathcal I_p[\rho]-\,d\,\mathcal E_p[\rho]\Big)_{|t=0}>0\,.
\]
\end{prop}
\noindent The function $\rho_0$ is explicitly constructed in~\cite{DEL-tlse}.

\smallskip To overcome the limitation $p\le2^\#$, one can consider a nonlinear diffusion of fast diffusion or porous medium type
\be{Eqn:FDE}
\frac{\partial\rho}{\partial t}=\Delta\rho^m\,.
\ee
With this flow, we no longer have $\frac d{dt}\mathcal E_p[\rho]=-\,\mathcal I_p[\rho]$ but we can still prove that
\[
\frac d{dt}\Big(\mathcal I_p[\rho]-\,d\,\mathcal E_p[\rho]\Big)\le0\,,
\]
for any $p\in[1,2^*]$. Proofs of the latter have been given in~\cite{MR2381156,1302}. We also refer to~\cite{DEKL2012,DEKL2014} for results which are more specific to the case of the sphere, and further references therein. Except for $p=1$ and $p=2^*$ with $d\ge3$, there is some flexibility in the choice of $m$, which can be used to build deficit functionals and improved inequalities: see~\cite{MR2381156,DEKL2014}. Notice that $\rho_0$ in Proposition~\ref{Prop:Counter-Example} is a function related with the nonlinear diffusion equation~\eqref{Eqn:FDE}.

The case of $\S^d$ highlights the limitations of linear flows and shows the flexibility and strength of nonlinear flows. At least for $p<2^*$, the optimal constant in~\eqref{Ineq:GNS} and~\eqref{Ineq:logSob} is established by proving that the minimum of $\mathcal I_p[\rho]-\,d\,\mathcal E_p[\rho]$ is $0$. Earlier results in~\cite{BV-V,MR1412446,MR1230930} can be reinterpreted as a purely elliptic method, which goes as follows. A positive minimizer actually exists by standard compactness arguments and any solution $\rho$ satisfies an Euler-Lagrange equation. By testing the equation with $\Delta\rho^m$, we observe that the solution is a constant and, as a consequence, that $\rho\equiv1$ because of the normalization. We will rely on a similar observation in the next two sections and refer to this method as the \emph{elliptic method.}

The method applies not only to minimizers, but also to any positive solution of the Euler-Lagrange equations. What we prove is a uniqueness result. Since constant functions are solutions, this proves that there are no non-constant solutions. This is why it is called a \emph{rigidity} result. 

Compared to~\cite{BV-V,MR1412446,MR1230930}, our approach provides a unified framework for \hbox{$p>2$} and $p<2$ (which is not covered in the above mentioned results). However, the main advantage of the method is that it explains why a local result (the best constant is given by the linearization around the constant functions) is actually global: $\mathcal I_p[\rho]-\,d\,\mathcal E_p[\rho]$ is strictly monotone decreasing under the action of the flow, unless the solution has reached the unique, trivial stationary state.

\AMSsubsection{Inequalities on compact manifolds}\labelSec{Sec:CompactManifolds}
The nonlinear diffusion flow method applies not only to spheres, but also to general compact manifolds. Without entering in the details, let us state a result of~\cite{1302}. Earlier important references are:~\cite{MR615628,BV-V,MR1631581,MR2381156}, among many other contributions which are listed in~\cite{1302}.

Let us assume that $(\mathcal M,g)$ is a smooth compact connected Riemannian manifold of dimension $d\ge1$, without boundary. We denote by $dv_g$ the volume element, by $\Lap$ the Laplace-Beltrami operator on $\mathcal M$, by~$\mathrm{Ric}$ the Ricci tensor and assume for simplicity that $\mathrm{vol}_g(\mathcal M)=1$. Let~$\lambda_1$ be the lowest positive eigenvalue of $-\Lap$ and
\[
\lambda_\star:=\inf_{u\in\mathrm H^2\,(\mathcal M)}\kern -4pt\frac{\displaystyle\iM{\Big[(1-\theta)\,(\Lap u)^2+\tfrac{\theta\,d}{d-1}\,\mathrm{Ric}(\nabla u,\nabla u)\Big]}}{\iM{|\nabla u|^2}}\,,\quad\theta=\frac{(d-1)^2\,(p-1) }{d\,(d+2)+p-1}\,.
\]
\begin{thm}\label{Thm:MainManifold} With the above notations, if $0<\lambda<\lambda_\star$, 
then for any $p\in(1,2)\cup(2,2^*)$, the equation
\[
-\,\Lap v+\frac\lambda{p-2}\,\(v-v^{p-1}\)=0
\]
has a unique positive solution in $C^2(\mathcal M)$, which is constant and equal to $1$.
\end{thm}
It has been shown in~\cite{1302} that nonlinear diffusion flows provide a unified framework for \emph{elliptic rigidity} and \emph{carr\'e du champ} methods. The computations heavily rely on the Bochner-Lichnerowicz-Weitzenb\"ock formula
\[
\tfrac12\,\Delta\,(|\nabla f|^2)=\|\mathrm {Hess}f\|^2+\nabla\cdot(\Delta f)\cdot\nabla f+\mathrm{Ric}(\nabla f, \nabla f)\,.
\]
More general results can be established using the so-called $CD(\rho,N)$ condition (see \cite{MR3155209} and references therein), but they are formal in most of the cases covered only by nonlinear flows. In dimension $d=2$, the Moser-Trudinger-Onofri inequality replaces in a certain sense Sobolev's inequality, and it is possible to extend the method described above to cover this case: see~\cite{Dolbeault20173059}. Bounded convex domains in $\R^d$ have also been considered in~\cite{1407} in relation with the Lin-Ni conjecture (homogeneous Neumann boundary conditions). Concerning unbounded domains, subcritical Ga\-gliardo-Niren\-berg have been established in the case of the line in~\cite{Dolbeault06082014} while R\'enyi entropy powers, which will be essential in Section~\refSec{Sec:DELM}, can be used in $\R^d$ to get sharp interpolation inequalities: see~\cite{MR3200617,Toscani-2014,DT-2015}.

\AMSsection{Rigidity on cylinders and sharp symmetry results in critical Caffarelli-Kohn-Nirenberg inequalities}\labelSec{Sec:CKN}

In this section we use a nonlinear flow to prove rigidity results for nonlinear elliptic problems on non-compact manifolds: cylinders and weigthed Euclidean spaces. All results of this section, and their proofs, can be found in~\cite{DEL-2015}.

\AMSsubsection{Three equivalent rigidity results}\labelSec{sec:1}
Let us consider the spherical cylinder $\mathcal C:=\mathbb R\times\S^{d-1}$ and denote by $s\in\R$ and $\omega\in\S^{d-1}$ the coordinates. Let $\Delta_\omega$ denote the Laplace-Beltrami operator on $\S^{d-1}$.
\begin{thm}\label{theorem:LV} Let $d\ge2$. For all $p\in(2,2^*)$ and $0<\Lambda\le\Lambda_{\rm FS}:=4\,\frac{d-1}{p^2-4}$, any positive solution $\varphi\in\mathrm H^1(\mathcal C)$ of
\be{eqlinder}
-\,\partial^2_s\,\varphi-\,\Delta_\omega\,\varphi+\Lambda\,\varphi=\varphi^{p-1}\quad\mbox{in}\quad\mathcal C
\ee
is, up to a translation in the $s$-direction, equal to
\[\label{varphistar}
\varphi_\Lambda(s):=\left(\tfrac p2\,\Lambda\right)^\frac1{p-2}\,\left(\cosh\left(\tfrac{p-2}2\,\sqrt\Lambda\,s\right)\right)^{-\frac2{p-2}}\quad\forall\,s\in\R\,.
\]
For any $\Lambda>\Lambda_{\rm FS}$, there are also positive solutions which do not depend only on $s$.
\end{thm}
A similar rigidity result holds for non-spherical cylinders $\mathbb R\times\mathfrak M$ where $\mathfrak M$ is a compact manifold, but in this case we cannot characterize the optimal set of parameters $\Lambda$ with our method: see~\cite{DEL-2015}.

Let 
\[\label{FS}
a_c:=\frac{d-2}2\quad\mbox{and}\quad b_{\rm FS}(a):=\frac{d\,(a_c-a)}{2\,\sqrt{(a_c-a)^2+d-1}}+a-a_c\,.
\]
By using the Emden-Fowler transformation
\be{EF}
v(r,\omega)=r^{a-a_c}\,\varphi(s,\omega)\quad\mbox{with}\quad r=|x|\,,\quad s=-\log r\quad\mbox{and}\quad\omega=\frac xr\,,
\ee
Theorem~\ref{theorem:LV} is equivalent to the following result.
\begin{thm}\label{theorem:Rigidity} Assume that $d\ge2$, $a<a_c$ and $\min\{a,b_{\rm FS}(a)\}<b\le a+1$. Then any nonnegative solution $v$ of
\be{EL}
-\,\nabla\cdot\big(|x|^{-2\,a}\,\nabla v\big)=|x|^{-b\,p}\,|v|^{p-2}\,v\,\quad\mbox{in}\quad\mathbb R^d\setminus\{0\}
\ee
which satisfies $\int_{\R^d}{\frac{|v|^p}{|x|^{b\,p}}\\,dx}<\infty$, is, up to a scaling, equal to
\[
v_\star(x)=\left(1+|x|^{(p-2)\,(a_c-a)}\right)^{-\frac2{p-2}}\quad\forall\,x\in\mathbb R^d\,.
\]
If $a<0$ and $a<b<b_{\rm FS}(a)$, there are also positive solutions which do not depend only on~$|x|$.\end{thm}

Let us define $\alpha_{\rm FS}:=\sqrt{\frac{d-1}{n-1}}$ and pick $n$ and $\alpha$ such that
\[
n=\frac{d-b\,p}\alpha=\frac{d-2\,a-2}\alpha+2=\frac{2\,p}{p-2}\,,
\]
so that we also have $p=2\,n/(n-2)$. Next we consider the diffusion operator
\[
\mathcal L\,w:=\alpha^2\(w''+\frac{n-1}r\,w'\)-\frac1{r^2}\,\Delta_\omega\,w\,.
\]
Then, with the change of variables
\[\label{wv}
v(r,\omega)=w(r^\alpha,\omega)\quad\forall\,(r,\omega)\in\mathbb R^+\times\S^{d-1}\,,
\]
Theorem~\ref{theorem:Rigidity} is equivalent to
\begin{thm}\label{theorem:Rigidity-dimension-n} Assume that $n>d\ge2$ and $p=2\,n/(n-2)$. If $0<\alpha\le\alpha_{\rm FS}$, then any nonnegative solution $w(x)=w(r,\omega)$ with $r\in\mathbb R_+$ and $\omega\in\S^{d-1}$ of
\be{EL-dimension-n}
-\,\mathcal L\,w=w^{p-1}\quad\mbox{in}\quad\mathbb R^d\setminus\{0\}
\ee
which satisfies $\int_{\R^d}{|x|^{n-d}\,|w|^p\,dx}<\infty$, is equal, up to a scaling, to
\[
w_\star(x)=\left(1+|x|^2\right)^{-\frac{n-2}2}\quad\forall\,x\in\mathbb R^d\,.
\]
If $\alpha>\alpha_{\rm FS}$, there are also solutions which do not depend only on $|x|$.\end{thm}
\noindent Let us complement these results with some remarks:\\
(i) If $n$ is an integer, then~\eqref{EL-dimension-n} is the Euler-Lagrange equation associated with the standard Sobolev inequality
\[
-\,\alpha^2\,\Delta w=w^\frac{n+2}{n-2}\quad\mbox{in}\quad\R^n\,,
\]
where $\Delta$ denotes the Laplacian operator in $\R^n$, but in the class of functions which depend only on the first $d-1$ angular variables.\\
(ii) The conditions on the parameters in Theorems~\ref{theorem:LV},~\ref{theorem:Rigidity} and~\ref{theorem:Rigidity-dimension-n} are equivalent:
\[
0<\Lambda\le\Lambda_{\rm FS}\;\Longleftrightarrow\;b_{\rm FS}^{-1}(b)\le a<a_c\;\Longleftrightarrow\;0<\alpha\le\alpha_{\rm FS}\,.
\]
(iii) Solutions of~\eqref{varphistar},~\eqref{EL} and~\eqref{EL-dimension-n} are stable (in a sense defined below) among non-symmetric solutions, \emph{i.e.}, solutions which explicitly depend on $\omega$, if and only if the above condition on the parameters is satisfied. Such a condition has been introduced in~\cite{Catrina-Wang-01}, but the sharp condition was established by V.~Felli and M.~Schneider in~\cite{Felli-Schneider-03}, and this is why we use the notation $\Lambda_{\rm FS}$, $b_{\rm FS}$ and $\alpha_{\rm FS}$ (see Section~\ref{BifurcationCritical}). Notice that stability is a local property while our uniqueness (rigidity) results are global.

\AMSsubsection{Optimal symmetry range in critical Caffarelli-Kohn-Nirenberg inequalities}\labelSec{sec:2}

The Caffarelli-Kohn-Nirenberg inequalities
\be{CKN}
\left(\int_{\R^d}{\frac{|v|^p}{|x|^{b\,p}}\,dx}\right)^{2/p}\le\,\mathsf C_{a,b}\int_{\R^d}{\frac{|\nabla v|^2}{|x|^{2\,a}}\,dx}\quad\forall\,v\in\mathcal D_{a,b}
\ee
appear in~\cite{Caffarelli-Kohn-Nirenberg-84}, under the conditions that $a\le b\le a+1$ if $d\ge3$, $a<b\le a+1$ if $d=2$, $a+1/2<b\le a+1$ if $d=1$, and $a<a_c$ where the exponent
\[\label{exponent-relationabp}
p=\frac{2\,d}{d-2+2\,(b-a)}
\]
is determined by the invariance of the inequality under scalings. Here $\mathsf C_{a,b}$ denotes the optimal constant in~\eqref{CKN} and the space $\mathcal D_{a,b}$ is defined by
\[
\mathcal D_{a,b}:=\Big\{\,v\in\mathrm L^p\big(\R^d,|x|^{-b}\,dx\big)\,:\,|x|^{-a}\,|\nabla v|\in\mathrm L^2\big(\R^d,dx\big)\Big\}\,.
\]
These inequalities were apparently introduced first by V.P. Il'in in~\cite{Ilyin} but are more known as \emph{Caffarelli-Kohn-Nirenberg inequalities}, according to~\cite{Caffarelli-Kohn-Nirenberg-84}. Up to a scaling and a multiplication by a constant, any extremal function for the above inequality is a nonnegative solution of~\eqref{EL}. It is therefore natural to ask whether $v_\star$ realizes the equality case in~\eqref{CKN}. Let
\[
\mathsf C_{a,b}^\star:=\frac{\left(\int_{\R^d}{\frac{|v_\star|^p}{|x|^{b\,p}}\,dx}\right)^{2/p}}{\int_{\R^d}{\frac{|\nabla v_\star|^2}{|x|^{2\,a}}\,dx}}=\tfrac p2\,|\S^{d-1}|^{1-\frac2p}\,(a-a_c)^{1+\frac2p}\(\frac{2\,\sqrt\pi\;\Gamma\big(\frac p{p-2}\big)}{(p-2)\,\Gamma\big(\frac{3\,p-2}{2\,(p-2)}\big)}\)^\frac{p-2}p\,.
\]
It was proved in~\cite{Felli-Schneider-03} that whenever $a<0$ and $b<b_{\rm FS}(a)$, the solutions of~\eqref{EL} are not radially symmetric: this is a \emph{symmetry breaking} result, based on the linear instability of $\mathcal F[v]:=\mathsf C_{a,b}^\star\int_{\R^d}{\frac{|\nabla v|^2}{|x|^{2\,a}}\,dx}-\big(\int_{\R^d}{\frac{|v|^p}{|x|^{b\,p}}\,dx}\big)^{2/p}$ at $v=v_\star$. The main \emph{symmetry} result of~\cite{DEL-2015} is
\begin{cor} Assume that $d\ge2$, $a<a_c$, and $b_{\rm FS}(a)\le b\le a+1$ if $a<0$. Then $\mathsf C_{a,b}=\mathsf C_{a,b}^\star$ and equality in~\eqref{CKN} is achieved by a function $v\in\mathcal D_{a,b}$ if and only if, up to a scaling and a multiplication by a constant, $v=v_\star$.\end{cor}
In other words, whenever $\mathcal F[v]$ is linearly stable at $v=v_\star$, then $v_\star$ is a global extremal function for~\eqref{CKN}.

\AMSsubsection{Sketch of the proof of Theorem~\ref{theorem:Rigidity-dimension-n}}\labelSec{sec:3}
The case $d=2$ requires some specific estimates so we shall assume that $d\ge 3$ for simplicity. Let
\be{uw}
u^{\frac 12-\frac1n}=w\quad\Longleftrightarrow\quad u=w^p\quad\mbox{with}\quad p=\frac{2\,n}{n-2}\,.
\ee

Up to a multiplicative constant, the right hand side in~\eqref{CKN} is transformed into a generalized \emph{Fisher information} functional
\be{Fisher}
\mathcal I[u]:=\irdmu{u\,|{\mathsf D}\mathsf p|^2}\quad\mbox{where}\quad\mathsf p=\frac m{1-m}\,u^{m-1}\,.
\ee
Here $d\mu=|x|^{n-d}\,dx$, $\mathsf p$ is the \emph{pressure function}, $\mathsf{D\,p}:=\big(\alpha\,\frac{\partial \mathsf p}{\partial r},\frac1r\,\nabla_\omega \mathsf p\big)$, and $\mathsf p'=\frac{\partial\mathsf p}{\partial r}$ and $\nabla_\omega \mathsf p$ respectively denote the radial and the angular derivatives of $\mathsf p$. The left hand side in~\eqref{CKN} is now proportional to a \emph{mass} integral, $\irdmu u$. In this section we consider the \emph{critical case} and make the choice $m=1-1/n$.

After these preliminaries, let us introduce the \emph{fast diffusion} flow
\be{FDE}
\frac{\partial u}{\partial t}=\L u^m\,,\quad m=1-\frac1n\,,
\ee
where the operator $\L$, which has been considered in Theorem~\ref{theorem:Rigidity-dimension-n}, is such that $\L w:=-\,{\mathsf D}^*\,{\mathsf D}\,w$. The flow associated with~\eqref{FDE} preserves the mass. At formal level, the key idea is to prove that $\mathcal I[u(t,\cdot)]$ is decreasing w.r.t.~$t$ if $u$ solves~\eqref{FDE}, and that the limit is $\mathcal I[w_\star^p]$. A long computation indeed shows that, if $u$ is a smooth solution of~\eqref{FDE} with the appropriate behavior as $x\to0$ and as $|x|\to+\infty$, then
\[
\frac d{dt}\mathcal I[u(t,\cdot)]\le-\,2\irdmu{\mathsf K[\mathsf p(t,\cdot)]\,u(t,\cdot)^m}
\]
where, with $r=|x|$, we have
\begin{multline}\label{Eqn:K}
\mathsf K[\mathsf p]=\alpha^4\left(1-\frac1n\right)\left[\mathsf p''-\frac{\mathsf p'}r-\frac{\Delta_\omega\,\mathsf p}{\alpha^2\,(n-1)\,r^2}\right]^2+2\,\alpha^2\,\frac1{r^2}\left|\nabla_\omega\mathsf p'-\frac{\nabla_\omega\mathsf p}r\right|^2\\
+(n-2)\big(\alpha_{\rm FS}^2-\alpha^2\big)\,\frac{|\nabla_\omega\mathsf p|^2}{r^4}+\zeta_\star\,(n-d)\,\frac{|\nabla_\omega\mathsf p|^4}{r^4}
\end{multline}
for some positive constant $\zeta_\star$. Hence, if $\alpha\le\alpha_{\rm FS}$, then $\mathcal I[u(t,\cdot)]$ is nonincreasing along the flow of~\eqref{FDE}. However, regularity and decay estimates needed to justify such computations are not known yet and this parabolic approach is therefore formal. As in Section~\refSec{Sec:Flows}, we can instead rely on an elliptic method, which can be justified as follows.

If $u_0$ is a nonnegative critical point of $\mathcal I$ under mass constraint, then
\[
0=\mathcal I'[u_0]\cdot \L u_0^m=\frac{d\mathcal I[u(t,\cdot)]}{dt}_{| t=0}\le-\,2\irdmu{\mathsf K[\mathsf p_0]\,\mathsf u_0^{1-n}}
\]
if $u$ solves~\eqref{FDE} with initial datum $u_0$. Here $\mathcal I'[u_0]$ denotes the differential of~$\mathcal I$ at~$u_0$. With $\mathsf p_0=\mathsf p(0,\cdot)$, this proves that $\nabla_\omega\mathsf p_0=0$: $\mathsf p_0$ is radially symmetric. By solving $\mathsf p_0''-\mathsf p_0'/r=0$, we obtain that $\mathsf p_0(x)=a+b\,|x|^2$ for some constants $a$, \hbox{$b\in\R^+$}. The conclusion easily follows.
\begin{prop}\label{Prop:DerivFisher} Let $w$ be a nonnegative solution of~\eqref{EL-dimension-n} and $\mathsf p=(n-1)\,w^{-\frac2{n-2}}$. Under the assumptions of Theorem~\ref{theorem:Rigidity-dimension-n}, if $\alpha\le\alpha_{\rm FS}$, then $\mathsf K[\mathsf p]=0$.\end{prop}
In practice, we prove that any solution of~\eqref{eqlinder} on $\mathcal C$ has good decay properties as $s\to\pm\infty$, by delicate elliptic estimates, which rely on the fact that $p=2\,n/(n-2)<2\,d/(d-2)$ is a subcritical exponent on the $d$-dimensional manifold $\mathcal C$. This is enough to justify all integrations by parts and prove as a consequence that a nonnegative solution of~\eqref{EL-dimension-n} satisfies $\mathsf K[\mathsf p]=0$: the conclusion follows as above. Notice that this amounts to test \eqref{EL-dimension-n} by $\L w^{2\,(n-1)/(n-2)}$.

\AMSsection{Rigidity and sharp symmetry results in subcritical Caffarelli-Kohn-Nirenberg inequalities}\labelSec{Sec:DELM}

In this section we consider a class of subcritical Caffarelli-Khon-Nirenberg inequalities and extend the results obtained for the critical case. Most results of this section have been published in~\cite{DELM}, a joint paper of the authors with M. Muratori.

\AMSsubsection{Subcritical Caffarelli-Kohn-Nirenberg inequalities}\labelSec{Sec:DELM-CKN}
With the notation
\[
\nrmRd w{q,\gamma}:=\(\ird{|w|^q\,|x|^{-\gamma}}\)^{1/q}\,,\quad\nrmRd wq:=\nrmRd w{q,0}\,,
\]
we define $\mathrm L^{q,\gamma}(\R^d)$ as the space $\{w\in\mathrm L^1_{\rm{loc}}(\R^d\setminus\{0\})\,:\,\nrmRd w{q,\gamma}<\infty\}$. We shall work in the space $\mathrm H^p_{\beta,\gamma}(\R^d)$ of functions $w\in\mathrm L^{p+1,\gamma}(\R^d)$ such that $\nabla w\in\mathrm L^{2,\beta}(\R^d)$, which can also be defined as the completion of $\mathcal D(\R^d\setminus\{0\})$ with respect to the norm
\[
\|w\|^2:=(p_\star-p)\,\nrmRd w{p+1,\gamma}^2+\nrmRd{\nabla w}{2,\beta}^2\,.
\]
Let us consider the family of subcritical \emph{Caffarelli-Kohn-Nirenberg interpolation inequalities} that can be found in~\cite{Caffarelli-Kohn-Nirenberg-84} and which is given by
\be{CKN-DELM}
\nrmRd w{2p,\gamma}\le\mathcal C_{\beta,\gamma,p}\,\nrmRd{\nabla w}{2,\beta}^\vartheta\,\nrmRd w{p+1,\gamma}^{1-\vartheta}\quad\forall\,w\in\mathrm H^p_{\beta,\gamma}(\R^d)\,.
\ee
Here the parameters $\beta$, $\gamma$ and $p$ are subject to the restrictions
\be{parameters}
d\ge2\,,\quad\gamma-2<\beta<\frac{d-2}d\,\gamma\,,\quad\gamma\in(-\infty,d)\,,\quad p\in\(1,p_\star\right]
\ee
with
\[
p_\star:=\frac{d-\gamma}{d-\beta-2}\quad\mbox{and}\quad\vartheta=\frac{(d-\gamma)\,(p-1)}{p\,\big(d+\beta+2-2\,\gamma-p\,(d-\beta-2)\big)}\,.
\]
The critical case $p=p_\star$ determines \hbox{$\vartheta=1$} and has been dealt with in Section~\refSec{Sec:CKN}, so we shall focus on the subcritical case $p<p_\star$. Here by \emph{critical} we simply mean that $\nrmRd w{2p,\gamma}$ scales like $\nrmRd{\nabla w}{2,\beta}$ and $\mathcal C_{\beta,\gamma,p}$ denotes the optimal constant in~\eqref{CKN-DELM}. The limit case $\beta=\gamma-2$ and $p=1$, which is an endpoint for~\eqref{parameters}, corresponds to Hardy-type inequalities: optimality is achieved among radial functions but there is no extremal function: see~\cite{0902}. The other endpoint is $\beta=(d-2)\,\gamma/d$, in which case $p_\star=d/(d-2)$: according to~\cite{Catrina-Wang-01} (also see Section~\ref{BifurcationCritical}), either $\gamma\ge0$, symmetry holds and there exists a symmetric extremal function, or $\gamma<0$, and then symmetry is broken but there is no extremal function. in all other cases, the existence of an extremal function for~\eqref{CKN-DELM} follows from standard methods: see~\cite{Catrina-Wang-01,DE2010,DMN2015} for related results.

When $\beta=\gamma=0$,~\eqref{CKN-DELM} is a Gagliardo-Nirenberg interpolation inequality which is well known to be related to the fast diffusion equation $\frac{\partial u}{\partial t}=\Delta u^m$ in $\R^d$, not only for $m=1-1/d$ but also for any $m\in[1-1/d,1)$. Here we generalize this observation to the weighted spaces.

\smallskip\emph{Symmetry} in~\eqref{CKN-DELM} means that the equality case is achieved by Aubin-Talenti type functions
\[
w_\star(x)=\(1+|x|^{2+\beta-\gamma}\)^{-1/(p-1)}\quad\forall\,x\in\R^d\,.
\]
On the contrary, there is \emph{symmetry breaking} if this is not the case, because the equality case is then achieved by a non-radial extremal function. It has been proved in~\cite{2016arXiv160208319B} that \emph{symmetry breaking} holds in~\eqref{CKN-DELM}~if
\be{set-symm-breaking}
\gamma<0\quad\mbox{and}\quad\beta_{\rm FS}(\gamma)<\beta<\frac{d-2}d\,\gamma
\ee
where
\[
\beta_{\rm FS}(\gamma):=d-2-\sqrt{(\gamma-d)^2-4\,(d-1)}\,.
\]
Under Condition~\eqref{parameters}, \emph{symmetry} holds in the complement of the set defined by~\eqref{set-symm-breaking}.
\begin{thm}\label{Thm:Main}{\sl Assume that~\eqref{parameters} holds and that
\be{Symmetry condition}
\beta\le\beta_{\rm FS}(\gamma)\quad\mbox{if}\quad\gamma<0\,.
\ee
Then the extremal functions for~\eqref{CKN-DELM} are radially symmetric and, up to a scaling and a multiplication by a constant, equal to $w_\star$.}\end{thm}
This means that~\eqref{set-symm-breaking} is the sharp condition for \emph{symmetry breaking}.

\AMSsubsection{A rigidity result}\labelSec{Sec:DELM-rigidity}
Up to a scaling and a multiplication by a constant, the Euler-Lagrange equation
\be{ELeq}
-\,\mbox{div}\,\big(|x|^{-\beta}\,\nabla w\big)=|x|^{-\gamma}\,\big(w^{2p-1}-\,w^p\big)\quad\mbox{in}\quad\R^d\setminus\{0\}
\ee
is satisfied by any extremal function for~\eqref{CKN-DELM}. In the range of parameters given by~\eqref{parameters} and~\eqref{Symmetry condition}, our method establishes the symmetry of all positive solutions.
\begin{thm}\label{Thm:Rigidity}{\sl Assume that~\eqref{parameters} and~\eqref{Symmetry condition} hold. Then all positive solutions to~\eqref{ELeq} in $\mathrm H^p_{\beta,\gamma}(\R^d)$ are radially symmetric and, up to a scaling, equal to $w_\star$.}\end{thm}
This is again a \emph{rigidity} result. Nonnegative solutions to~\eqref{ELeq} are actually positive by the standard Strong Maximum principle. Theorem~\ref{Thm:Main} is therefore a consequence of Theorem~\ref{Thm:Rigidity}.

\AMSsubsection{Sketch of the proof of Theorem~\ref{Thm:Rigidity}}\labelSec{sec:TT3}
Let us give an outline of the strategy of~\cite{DELM}. As in the critical case, Inequality~\eqref{CKN-DELM} for a function $w$ can be transformed by the change of variables
\[
w(x)=v(r^\alpha,\omega)\,,
\]
where $r=|x|\neq0$ and $\omega=x/r$, in the new inequality 
\be{CKN-DELM2}
\(\irdmu{|v|^{2p}}\)^\frac1{2p}\le\mathcal K_{\alpha,n,p}\(\irdmu{|\DD v|^2}\)^\frac\vartheta2\(\irdmu{|v|^{p+1}}\)^\frac{1-\vartheta}{p+1}
\ee
with $\mathcal K_{\alpha,n,p}=\alpha^{-\zeta}\,\mathcal C_{\beta,\gamma,p}$, $\zeta=\frac\vartheta2+\frac{1-\vartheta}{p+1}-\frac1{2\,p}$ and $d\mu=|x|^{n-d}\,dx$. The condition for the change of variables is
\[
n=\frac{d-\beta-2}\alpha+2=\frac{d-\gamma}\alpha\,,
\]
which reflects the fact that the weights are all the same in~\eqref{CKN-DELM2}. It is solved by
\[\label{Eqn:alpha-n}
\alpha=1+\frac{\beta-\gamma}2\quad\mbox{and}\quad n=2\,\frac{d-\gamma}{\beta+2-\gamma}\,.
\]
Inequality~\eqref{CKN-DELM2} is a Caffarelli-Kohn-Nirenberg inequality with weight $|x|^{n-d}$ in all terms, and $\DD v:=\big(\alpha\,\frac{\partial v}{\partial s},\frac1s\,\nabla_{\kern-2pt\omega}v\big)$. Notice that $p_\star=\frac n{n-2}$, so that $2\,p_\star$ is the critical Sobolev exponent associated with the \emph{fractional dimension} $n$ considered in~\eqref{uw}.

With a generalized \emph{Fisher information}~$\mathcal I$ and the \emph{pressure function} $\mathsf p$ defined by~\eqref{Fisher}, we consider the \emph{subcritical range} $m_1:=1-1/n<m<1$. If $u$ is smooth solution of~\eqref{FDE} with sufficient decay properties, we obtain that $\mathcal I$ evolves according to
\[\label{BLWevol}
\frac d{dt}\mathcal I[u(t,\cdot)]=-\,2\irdmu{\mathcal R[\mathsf p(t,\cdot)]\,u(t,\cdot)^m}\quad\mbox{with}\quad\mathcal R[\mathsf p]:=\mathsf K[\mathsf p]+\(m-m_1\)\(\L\mathsf p\)^2\,,
\]
where $\mathsf K$ is given by~\eqref{Eqn:K}. We recover the result of the critical case of Section~\refSec{sec:3} by taking the limit as $m\to m_1$. 

Inspired by tools of Information Theory and~\cite{MR3200617,Toscani-2014,DT-2015}, we introduce the generalized \emph{R\'enyi entropy power} functional
\[
\mathcal F[u]:=\(\irdmu{u^m}\)^\sigma\quad\mbox{with}\quad\sigma=\frac2n\,\frac1{1-m}-1>1
\]
and observe that $\mathcal F''$ has the sign of $-\,\mathcal H[u(t,\cdot)]$ where
\[
\mathcal H[u]:=(m-m_1)\irdmu{\left|\L\mathsf p-\frac{\irdmu{u\,|\DD\mathsf p|^2\,u^m}}{\irdmu{u^m}}\right|^2}+\irdmu{\mathcal R[\mathsf p]\,u^m}\,.
\]
Here $\mathcal F'$ denotes the derivative with respect to $t$ of $\mathcal F[u(t,\cdot)]$. The computation requires many integrations by parts. The fact that boundary terms do not contribute can be justified if $u$ is a nonnegative critical point, \emph{i.e.}, a minimizer of $\mathcal F'$ under mass constraint. Indeed, the minimization of
\[
\(\irdmu{v^{p+1}}\)^{\sigma-1}\irdmu{|\DD v|^2}\quad\mbox{with}\quad v=u^{m-1/2}
\]
under the constraint that $\irdmu u=\irdmu{v^{2p}}$ takes a given positive value is equivalent to the \emph{Caffarelli-Kohn-Nirenberg interpolation inequalities}~\eqref{CKN-DELM}.

To make the argument rigorous, we can argue as in Section~\refSec{sec:3} by taking $u$ as initial datum and performing the computation of $\mathcal F''$ at $t=0$ only. In other words, we are simply testing the Euler-Lagrange equation satisfied by $u$ with $\L u^m$. By elliptic regularity (the estimates are as delicate as in the critical case and we refer to~\cite{DELM} for details), we have enough estimates to prove that $\mathcal H[u]=0$ and deduce that $\mathsf p(x)=\mathsf a+\mathsf b\,|x|^2$ for some real constants $\mathsf a$ and~$\mathsf b$.

\AMSsubsection{Considerations on the optimality of the method}
The \emph{symmetry breaking} condition in~\eqref{CKN} and~\eqref{CKN-DELM} has been established by proving the linear instability of radial critical points, in~\cite{Felli-Schneider-03} and~\cite{2016arXiv160208319B} respectively. This amounts to a spectral gap condition in a Hardy-Poincar\'e inequality: see~\cite{2016arXiv160208319B} for details. It is remarkable that the symmetry holds whenever radial critical points are linearly stable and this deserves an explanation. The solution of~\eqref{FDE} is attracted by self-similar Barenblatt functions as $t\to+\infty$. Since these Barenblatt functions are precisely the radial critical points of our variational problem, the asymptotic rate of convergence is determined by the previous spectral gap, in self-similar variables. It can be checked that the condition that appears in the \emph{carr\'e du champ} method, which amounts to prove that a quadratic form has a sign, is the same in the asymptotic regime as $t\to+\infty$ as the quadratic form which is used to check symmetry breaking. Hence either symmetry breaking occurs, or the \emph{carr\'e du champ} method shows that the \emph{R\'enyi entropy power} functional is monotone non-increasing, at least in the asymptotic regime: see~\cite{DEL-JEPE} for details. To conclude in the critical case, it is enough to observe that all terms in the expression of $\mathsf K[\mathsf p]$ in~\eqref{Eqn:K} are quadratic, except the last one, which has a sign and is negligible compared to the others in the asymptotic regime: the sign condition for $\mathsf K[\mathsf p]$ away from the asymptotic regime is the same as when $t\to+\infty$. This explains why our method for proving \emph{symmetry} gives the optimal range in the critical case. In the subcritical regime, a similar observation can also be done.

\AMSsection{Bifurcations and symmetry breaking}

The results of this section are taken mostly from~\cite{DETT,DE-num,DE-branches}.

\AMSsubsection{Rigidity and bifurcations}\labelSec{BifurcationCritical}
Let us come back to the critical Caffarelli-Kohn-Nirenberg inequality and consider the Emden-Fowler transformation~\eqref{EF}. As noted in~\cite{Catrina-Wang-01}, Inequality~\eqref{CKN} is transformed into the Gagliardo-Nirenberg-Sobo\-lev inequality
\[
\nrmC{\nabla\varphi}2^2+\Lambda\,\nrmC\varphi2^2\ge\mu(\Lambda)\,\nrmC\varphi p^2\quad\forall\,\varphi\in\mathrm H^1(\mathcal C)
\]
where $\mu(\Lambda)=\mathsf C_{a,b}^{-1}\,\big|\Sp^{d-1}\big|^{1-2/p}$. Here $\mathcal C:=\R\times\Sp^{d-1}$ is a cylinder and, as in Section~\refSec{Sec:sphere}, we adopt the convention that the measure on the sphere is the uniform probability measure. The extremal functions are, up to multiplication by a constant, and dilation, solutions of~\eqref{eqlinder}.

If we restrict the study to symmetric functions, that is, $v(r)=r^{a-a_c}\,\varphi(-\log r)$ with $r=|x|$, then the inequality degenerates into the simple Gagliardo-Niren\-berg-Sobo\-lev inequality
\[
\nrmR{\nabla\varphi}2^2+\Lambda\,\nrmR\varphi2^2\ge\mu_\star(\Lambda)\,\nrmR\varphi p^2\quad\forall\,\varphi\in\mathrm H^1(\R)\,.
\]
Here we denote by
\[
\mu_\star(\Lambda)=\mu_\star(1)\,\Lambda^\frac{p+2}{2\,p}
\]
the optimal constant and notice that $\varphi_\star(s)=\big(\frac12\,p\,\Lambda\,\cosh\big(\tfrac{p-2}2\,\sqrt{\Lambda}\,s\big)^{-2}\big)^{1/(p-2)}$ is an optimal function, which is the unique solution of $-\,\varphi''+\Lambda\,\varphi=|\varphi|^{p-2}\,\varphi$ on~$\R$, up to translations. With this notation, we have $\mu_\star(\Lambda)=\nrmR{\varphi_\star}p^{p-2}$. If we linearize
\[
\nrmC{\nabla\varphi}2^2+\Lambda\,\nrmC\varphi2^2-\mu_\star(\Lambda)\,\nrmC\varphi p^2
\]
around $\varphi=\varphi_\star$, V.~Felli and M.~Schneider found in~\cite{Felli-Schneider-03} that the lowest eigenvalue of the quadratic form, that is, the lowest positive eigenvalue of the P\"oschl-Teller operator $-\frac{d^2}{ds^2}+\Lambda+d-1-(p-1)\,\varphi_\star^{p-2}$, is given by $\lambda_1(\Lambda)=-\frac 14\,(p^2-4)\,(\Lambda-\Lambda_{\rm FS})$, so that $\lambda_1(\Lambda_{\rm FS})<0$ if and only if
\[
\Lambda>\Lambda_{\rm FS}:=4\,\frac{d-1}{p^2-4}\,.
\]
See~\cite[p.~74]{Landau-Lifschitz-67} for details. This condition is the \emph{symmetry breaking} condition of Theorem~\ref{theorem:LV}. The branch of non-radial solutions bifurcating from $\Lambda=\Lambda_{\rm FS}$ has been computed numerically in~\cite{DE-num} and an example is shown in Fig.~\ref{Fig1}. By construction, we know that $\Lambda\mapsto\mu(\Lambda)$ is increasing, concave, and we read from Theorem~\ref{theorem:LV} that the non-symmetric branch bifurcates from $\Lambda=\Lambda_{\rm FS}$, and is such that $\mu(\Lambda)<\mu_\star(\Lambda)$ if \hbox{$\Lambda>\Lambda_{\rm FS}$}. This simple scenario explains the symmetry and symmetry breaking properties in~\eqref{CKN}, but is not generic as we shall see next in the case of more complicated interpolation inequalities.
\setlength\unitlength{1cm}
\begin{figure}[!ht]\begin{center}\begin{picture}(7,5)
\put(0,0){\includegraphics[width=7cm]{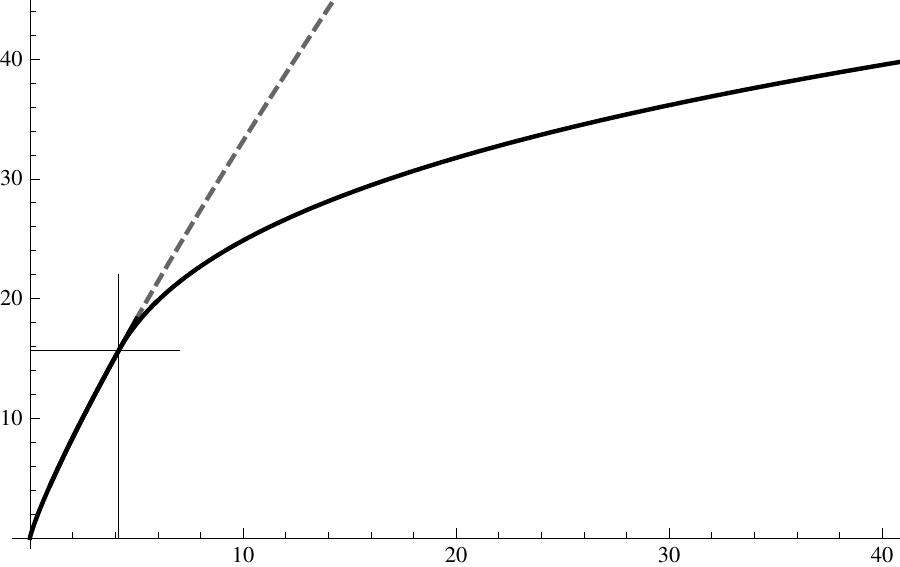}}
\put(6.7,0.45){$\Lambda$}
\put(6.4,3.4){$\mu(\Lambda)$}
\put(2.7,4.2){$\mu_\star(\Lambda)$}
\put(0.75,-0.15){$\Lambda_{\rm FS}$}
\put(-1.1,1.6){$\mu_\star(\Lambda_{\rm FS})$}
\end{picture}\caption{\label{Fig1}\small Branches for $p=2.8$, $d=5$, $\theta=1$.}
\end{center}\end{figure}
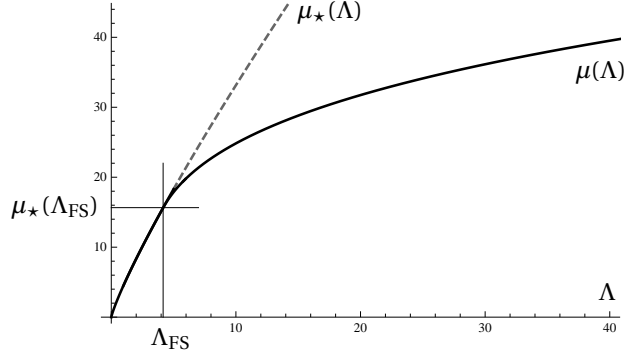

\AMSsubsection{Bifurcations, reparametrization and turning points}\labelSec{Sec:Bifurcations}
Let us consider the interpolation inequality
\be{Ineq:CKN}
\(\ird{\frac{|u|^p}{|x|^{bp}}}\)^\frac2p\le\mathsf C_{a,b,\theta}\(\ird{\frac{|\nabla u|^2}{|x|^{2a}}}\)^\theta\(\ird{\frac{|u|^2}{|x|^{2\,(a+1)}}}\)^{1-\theta}
\ee
with $d\ge1$, $p\in(2,2^*)$ or $p=2^*$ if $d\ge3$, and $\theta\in(\vartheta(p),1]$ with $\vartheta(p):=d\,\frac{p-2}{2\,p}$. The scaling invariance imposes $p=2\,d/\big(d-2+2\,(b-a)\big)$. As proved in~\cite{Caffarelli-Kohn-Nirenberg-84}, the above inequalities hold with a finite constant $\mathsf C_{a,b,\theta}$ if $a<a_c=(d-2)/2$, and $b\in(a+1/2,a+1]$ when $d=1$, $b\in(a,a+1]$ when $d=2$ and $b\in[a,a+1]$ when $d\ge3$. Moreover, there exist extremal functions for the inequalities~\eqref{Ineq:CKN} for any $p\in(2,2^*)$ and $\theta\in(\vartheta(p),1)$ or $\theta=\vartheta(p)$ and $d\ge2$, with $a_c-a>0$ not too large. On the contrary equality is never achieved for $p=2$, or $a<0$, $p=2^*$ and $d\ge 3$, or $d=1$ and $\theta=\vartheta(p,1)$. The existence of extremal functions has been studied in~\cite{DE2010}. We may notice that
\[
0\le\vartheta(p)\le\theta<1\quad\Longleftrightarrow\quad2\le p\le p^*(d,\theta):=\frac{2\,d}{d-2\,\theta}<2^*\,.
\]
With the same conventions as in the previous subsection, the Emden-Fowler change of variables~\eqref{EF} transforms~\eqref{Ineq:CKN} into the Gagliardo-Nirenberg-Sobolev inequality
\be{CKNthetaC}
\(\nrmC{\nabla\varphi}2^2+\Lambda\,\nrmC\varphi2^2\)^\theta\,\nrmC\varphi2^{2(1-\theta)}\ge\mu(\theta,\Lambda)\,\nrmC\varphi p^2\quad\forall\,\varphi\in\mathrm H^1(\mathcal C)
\ee
on $\mathcal C:=\R\times\Sp^{d-1}$, with $\Lambda=(a-a_c)^2$ and $\mu(\theta,\Lambda)=\mathsf C_{a,b,\theta}^{-1}\,\big|\Sp^{d-1}\big|^{1-2/p}$. Of course, the case $\theta=1$ corresponds to the critical case and, consistently, we write $\mu(1,\Lambda)=\mu(\Lambda)$.

For $\theta<1$, the Euler-Lagrange equation of an extremal function on $\mathcal C$ is
\be{EL2}
-\Delta\varphi+\frac1\theta\((1-\theta)\,\frac{\nrmcnd{\nabla\varphi}2^2}{\nrmcnd\varphi2^2}+\Lambda\)\varphi-\frac{\nrmcnd{\nabla\varphi}2^2+\Lambda\,\nrmcnd\varphi2^2}{\theta\,\nrmC\varphi p^p}\,\varphi^{p-1}=0\,.
\ee
Up to the reparametrization
\[
\Lambda\mapsto\lambda=\frac1\theta\Big[(1-\theta)\,t[\varphi]+\Lambda\Big]\quad\mbox{where}\quad t[\varphi]:=\frac{\nrmcnd{\nabla\varphi}2^2}{\nrmcnd\varphi2^2}
\]
and a multiplication by a constant, an extremal function $\varphi$ for~\eqref{CKNthetaC} solves~\eqref{eqlinder}. In other words, we can use the set of solutions in the critical case $\theta=1$ to parametrize the solutions corresponding to $\theta<1$.

Let us start with the symmetric functions. With an evident notation, we define $\mu_\star(\theta,\Lambda)$ as the optimal constant in the inequality corresponding to~\eqref{CKNthetaC} restricted to symmetric functions, \emph{i.e.}, functions depending only on $s\in\R$. If we denote by $\varphi_{\star,\lambda}$ the function
\[
\varphi_{\star,\lambda}(s)=\(\frac12\,p\,\lambda\,\cosh\big(\tfrac{p-2}2\,\sqrt{\lambda}\,s\big)^{-2}\)^\frac1{p-2}
\]
for any $\lambda>0$, then $t[\varphi_{\star,\lambda}]$ is explicit and we can parametrize the set $\big\{\big(\Lambda,\mu_\star(\theta,\Lambda)\big)\,:\,\Lambda>0\big\}$ by $\big\{\big(\theta\,\lambda-(1-\theta)\,t[\varphi_{\star,\lambda}],\mu_\star(\lambda)\big)\,:\,\lambda>0\big\}$. It turns out that the equation $\Lambda=\theta\,\lambda-(1-\theta)\,t[\varphi_{\star,\lambda}]$ can be inverted, which allows us to obtain $\lambda=\Lambda_*^\theta(\Lambda)$ and get an explicit expression for
\[
\mu_\star(\theta,\Lambda)=\mu_\star\big(\Lambda_*^\theta(\Lambda)\big)=\mu_\star\big(\Lambda_*^\theta(1)\big)\,\Lambda^{\theta-\frac{p-2}{2\,p}}\,.
\]
According to~\cite{DDFT}, a Taylor expansion around $\varphi_{\star,\Lambda_{\rm FS}}$ shows that for any $\Lambda>\Lambda_{\rm FS}^\theta$, where
\[
\Lambda_{\rm FS}^\theta:=\theta\,\mu_{\rm FS}-(1-\theta)\,t[u_{\rm FS}]\,,
\]
the function $\varphi_{\star,\lambda}$ with $\lambda=\Lambda_*^\theta(\Lambda)$ is linearly unstable, so that $\mu(\theta,\Lambda)<\mu_\star(\theta,\Lambda)$. 

\begin{figure}[!ht]\begin{center}\begin{picture}(7,5)
\put(-3,0){\includegraphics[width=8cm]{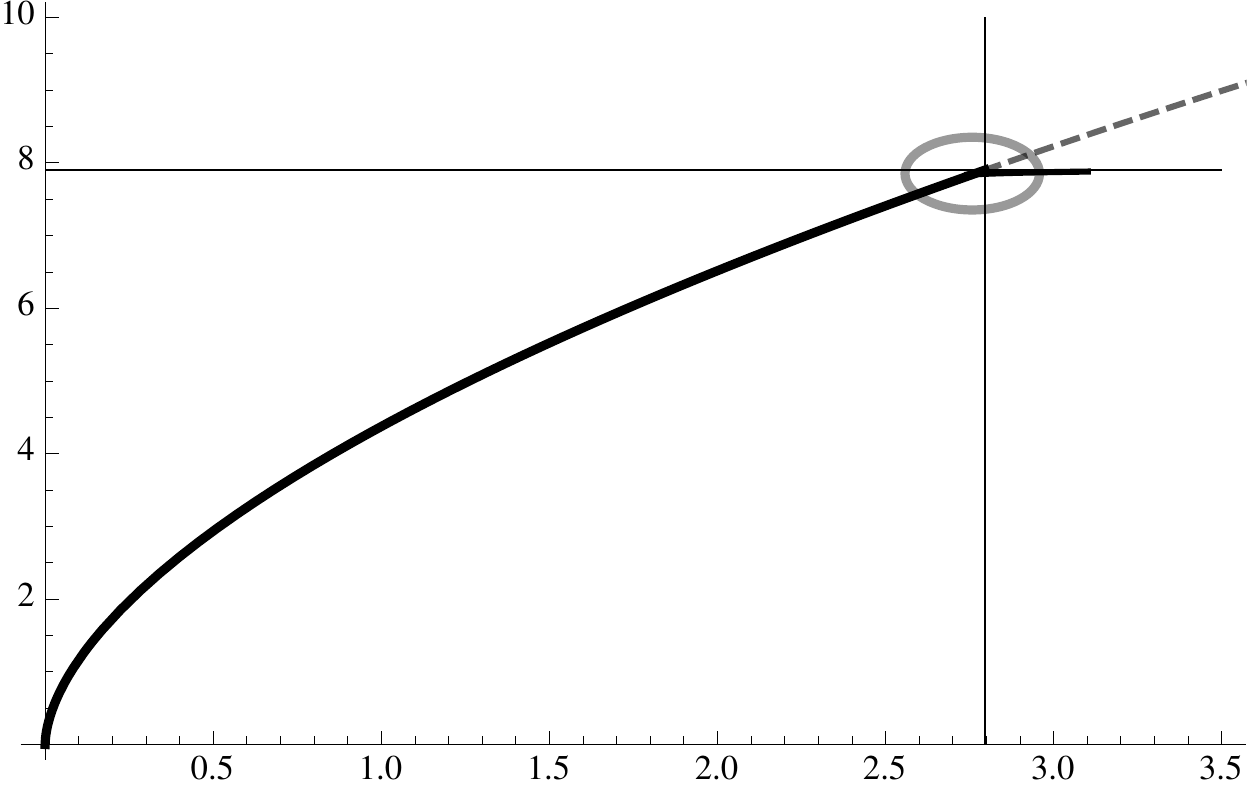}}
\put(5,0.5){\includegraphics[width=6cm]{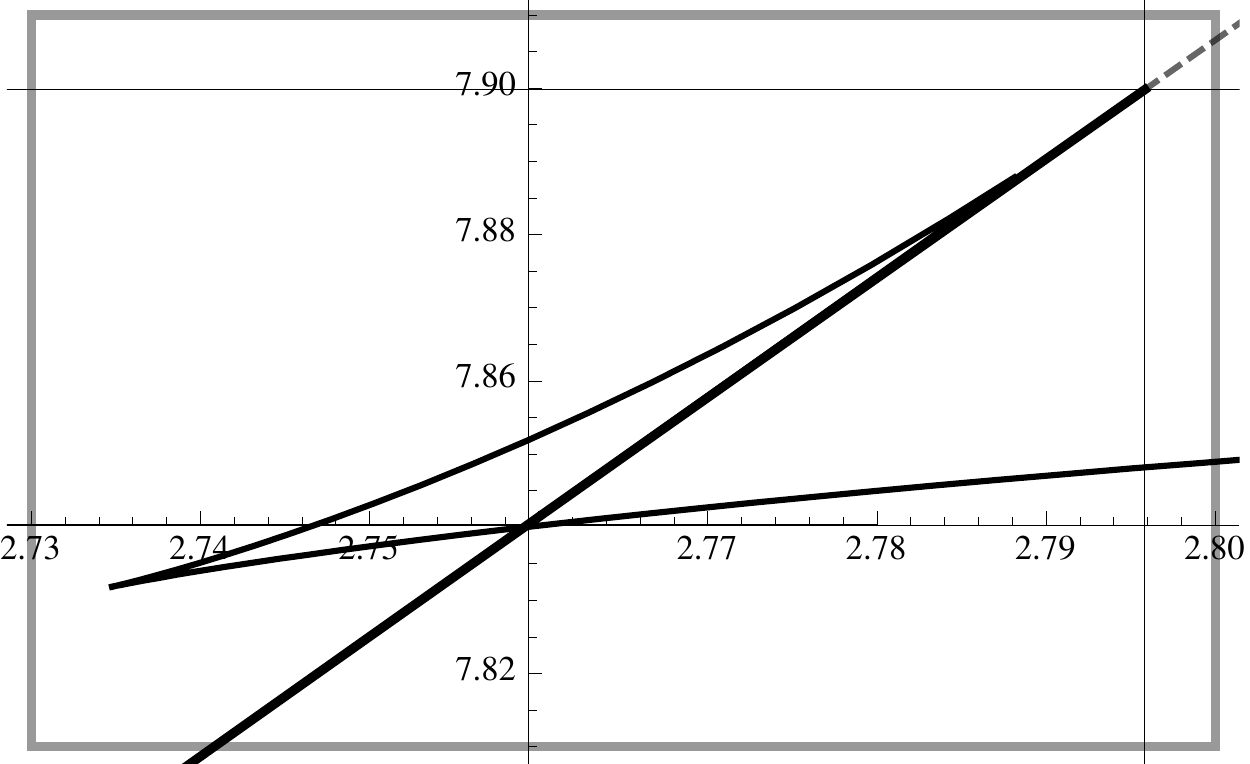}}
\put(4.5,0.5){$\Lambda$}
\put(3.7,3.5){$\mu(\theta,\Lambda)$}
\put(-0.5,3.1){$\mu_\star(\theta,\Lambda)$}
\put(3,-0.2){$\Lambda_{\rm FS}^\theta$}
\put(-4.5,3.7){$\mu_\star(\theta,\Lambda_{\rm FS}^\theta)$}
\put(6.8,2){$\courbe$}
\put(9,2){$\courbe$}
\put(6.5,1.2){$\courbe$}
\put(10.2,0.1){$\Lambda_{\rm FS}^\theta$}
\put(11,3.7){$\mu_\star(\theta,\Lambda_{\rm FS}^\theta)$}
\end{picture}
\vspace*{0.2cm}
\caption{\label{Fig2}\small Branches for $p=2.8$, $d=5$, $\theta=0.718$. Left: the bifurcation point $(\Lambda_{\rm FS}^\theta,\mu_\star(\theta,\Lambda_{\rm FS}^\theta)$ is at the intersection of the horizontal and vertical lines. The area enclosed in the small ellipse is enlarged in the right plot: the branch has a \emph{turning point} and $\mu(\theta,\Lambda_{\rm FS}^\theta)<\mu_\star(\theta,\Lambda_{\rm FS}^\theta)$.}
\end{center}\end{figure}
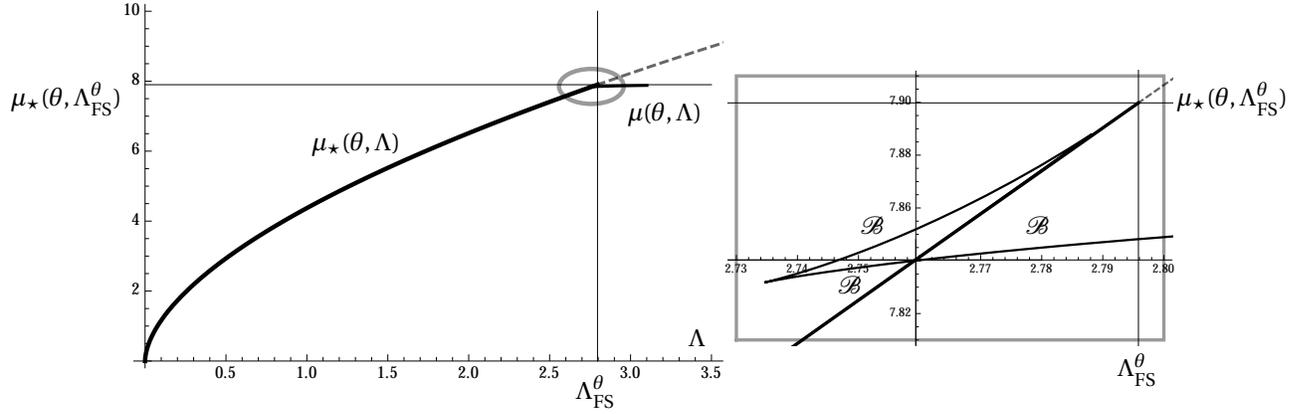

The case of non-symmetric functions is more subtle because we do not know the exact multiplicity of the solutions of~\eqref{eqlinder} in the symmetry breaking range. There is a branch of non-symmetric solutions of~\eqref{EL2} which bifurcates from the branch of symmetric solutions at $\Lambda=\Lambda_{\rm FS}^\theta$. This branch has been computed numerically in~\cite{DE-num} and a formal asymptotic expansion was performed in a neighborhood of the bifurcation point in~\cite{DE-branches}. Because of the reparametrization of the solutions of~\eqref{EL2} by the solutions of~\eqref{eqlinder}, we can use the branch $\lambda\mapsto\varphi_\lambda$ of non-symmetric extremal functions for $\lambda>\Lambda_{\rm FS}$ to get an upper bound of $\mu(\theta,\Lambda)$:
\[
\mu(\theta,\Lambda)\le\mu(\lambda)\quad\mbox{for any $\lambda>\Lambda_{\rm FS}$ such that}\quad\Lambda=\theta\,\lambda-(1-\theta)\,t[\varphi_\lambda]\,.
\]
Actually, we deduce from the branch $\lambda\mapsto\varphi_\lambda$ of non-symmetric extremal functions an entire branch of non-symmetric solutions of~\eqref{EL2} which is parametrized by $\lambda$ and deduce a parametric curve $\courbe:=\big\{\big(\Lambda(\lambda):=\theta\,\lambda-(1-\theta)\,t[\varphi_\lambda],\mu(\lambda)\big)\,:\,\lambda>\Lambda_{\rm FS}\big\}$ which can be used to bound $\mu(\theta,\Lambda)$ from above. If $(\Lambda,\mu)\in\courbe$, we have no proof that $\varphi_\lambda$ is optimal if $\mu(\lambda)<\mu_\star(\theta,\Lambda)$, but at least we know that
\[
\mu(\lambda)=\(\nrmC{\nabla\varphi_\lambda}2^2+\Lambda(\lambda)\,\nrmC{\varphi_\lambda}2^2\)^\theta\,\nrmC{\varphi_\lambda}2^{2(1-\theta)}\,\nrmC{\varphi_\lambda}p^{-2}\,.
\]
Some numerical results are shown in Fig.~\ref{Fig2}.

The formal asymptotic expansion of~\cite{DE-branches} suggests that there are only two possible generic scenarii:
\begin{enumerate}
\item[(i)] Either the curve $\courbe$ bifurcates to the right, that is, $\courbe$ is included in the region $\Lambda\ge\Lambda_{\rm FS}^\theta$, and $\Lambda\mapsto\mu(\theta,\Lambda)$ is qualitatively expected to be as in Fig.~\ref{Fig1}. We know that this is what happens for $\theta=1$ and expect a similar behavior for any $\theta$ close enough to $1$. In this case, the region of symmetry breaking is characterized by the linear instability of the symmetric optimal functions.
\item[(ii)] Or the curve $\courbe$ bifurcates to the left. For $\lambda-\Lambda_{\rm FS}>0$, small enough, the curve $\lambda\mapsto\big(\Lambda(\lambda),\mu(\lambda)\big)$ satisfies $\Lambda(\lambda)<\Lambda_{\rm FS}^\theta$ and $\mu(\lambda)>\mu_\star(\theta,\Lambda(\lambda))$. In that case, the region of symmetry breaking does not seem to be characterized by the linear instability of the symmetric optimal functions and we numerically observe a turning point as in Fig.~\ref{Fig2} (right).
\end{enumerate}
In~\cite{1406}, \emph{a priori} estimates for branches with $\theta<1$ were deduced from the known symmetry results (later improved in~\cite{DEL-2015}). This further constrains $\courbe$ and the symmetry breaking region and determines a lower bound for the value of $\Lambda$ corresponding to a turning point of the branch. There are many open questions concerning $\courbe$ and the set of extremal functions when $\theta<1$, but at least we can prove that the symmetry breaking range does not always coincide with the region of linear instability of symmetric optimal functions.

\AMSsubsection{Symmetry breaking and energy considerations}\labelSec{Sec:energy}
The exponent $\vartheta(p)$ is the exponent which appears in the Gagliardo-Nirenberg inequality
\be{GN}
\nrmRd{\nabla u}2^{2\,\vartheta(p)}\,\nrmRd u2^{2\,(1-\vartheta(p))}\ge\mathsf C_{\rm{GN}}(p)\,\nrmRd up^2\quad\forall\,u\in\mathrm H^1(\R^d)\,.
\ee
By considering an extremal function for this inequality and translations, for any $p\in(2,2^*)$, one can check that
\[
\mu(\vartheta(p),\Lambda)\le\mathsf C_{\rm{GN}}(p)\quad\forall\,\Lambda>0\,.
\]
\begin{lemma}\label{Lem:comparison} Let $d\ge2$. For any $p\in(2,2^*)$, if $\mathsf C_{\rm{GN}}(p)<\mu_\star(\vartheta(p),\Lambda_{\rm FS}^{\vartheta(p)})$, there exists $\Lambda_s\in(0,\Lambda_{\rm FS}^{\vartheta(p)})$ such that $\mu(\vartheta(p),\Lambda)=\mu_\star(\vartheta(p),\Lambda)$ if and only if $\Lambda\in(0,\Lambda_s]$.\end{lemma}
The fact that the symmetry range is an interval of the form $(0,\Lambda_s]$ can be deduced from a scaling argument: see~\cite{DELT09,DETT} for details. The result is otherwise straightforward but difficult to use because the value of $C_{\rm{GN}}(p)$ is not known explicitly. From a numerical point of view, it gives a simple criterion, which has been implemented in~\cite{DETT}. Moreover, in~\cite{DE-branches}, it has been observed numerically that the condition $\mathsf C_{\rm{GN}}(p)<\mu_\star(\vartheta(p),\Lambda_{\rm FS}^{\vartheta(p)})$ is equivalent to a \emph{bifurcation to the left} as in Fig.~\ref{Fig2}.

For $\theta$ and $p-2$ small enough, the assumption of Lemma~\ref{Lem:comparison} holds. Let us consider the Gaussian test function $\mathsf g(x):=(2\,\pi)^{-d/4}\,\exp(-|x|^2/4)$ in~\eqref{Ineq:CKN} and consider
\[
h(p):=\frac{\nrmRd{\nabla\mathsf g}2^{2\,\theta}\,\nrmRd{\mathsf g}2^{2\,(1-\theta)}}{\nrmRd{\mathsf g}p^2}\,\frac1{\mu_\star(\theta,\Lambda_{\rm FS}^\theta)}\quad\mbox{with}\quad\theta=\vartheta(p)\,.
\]
A computation shows that $\lim_{p\to2_+}h(p)=1$ and $\lim_{p\to2_+}\frac{dh}{dp}(p)<0$. For $p-2>0$, small enough, we obtain that
\[
\mathsf C_{\rm{GN}}(p)\le h(p)<\mu_\star(\theta,\Lambda_{\rm FS}^\theta)\,.
\]
A perturbation argument has been used in~\cite{DETT} to establish the following result.
\begin{thm}\label{Thm:DETT} Let $d\ge2$. There exists $\eta>0$ such that for any $p\in(2,2+\eta)$,
\[
\mu(\theta,\Lambda)<\mu_\star(\theta,\Lambda)\quad\mbox{if}\quad\Lambda_{\rm FS}^\theta-\eta<\Lambda<\Lambda_{\rm FS}^\theta\quad\mbox{and}\quad\vartheta(p)<\theta<\vartheta(p)+\eta\,.
\]
\end{thm}

\AMSsubsection{An open question}
The criterion considered in Lemma~\ref{Lem:comparison} is based on energy considerations and provides only a sufficient condition for symmetry breaking. It is difficult to check it in practice, except in asymptotic regimes of the parameters. The formal expansions of the branch near the bifurcation points are based on a purely local analysis, and suggest another criterion: either the branch bifurcates to the right and the symmetry breaking range is characterized by the linear instability of the symmetric optimal functions, or the branch bifurcates to the left, and this is not anymore the case. Is such an observation, which has been made numerically only for some specific values of~$p$, true in general? This seems to be true when $\theta$ is close enough to~$\vartheta(p)$ and at least in this regime we can conjecture that \emph{the symmetry breaking range is not characterized by the linear instability of the symmetric optimal functions if and only if the branch bifurcates to the left}.

An additional question, which corresponds to a limiting case, goes as follows. If $\theta=\vartheta(p)$, is the range of symmetry determined exactly by the value of the optimal constant in~\eqref{GN}, when it is below $\mu_\star(\theta,\Lambda_{\rm FS}^\theta)$? Numerically, this is supported by the fact that, in this case, the curve $\courbe$ is monotone increasing as a function of~$\Lambda$.
 
In the study of the symmetry issue in~\eqref{CKN} and~\eqref{CKN-DELM}, the key tool is the nonlinear flow, which extends a local result (linear stability) to a global result (rigidity). A similar tool would be needed to answer the conjecture. In the case $\theta=\vartheta(p)$, it would be crucial to obtain a variational characterization of the non-symmetric solutions in the curve of non-symmetric functions $\courbe$ and a uniqueness result for any given $\Lambda$.


\def\refname{}
\begin{center}\sc{References}\end{center}\vspace*{-1.5cm}
\bibliographystyle{siam}
\setstretch{0.8}
\small\begin{spacing}{0.9}
\bibliography{ICM}
\end{spacing}

\begin{center}\vspace*{-6pt}\rule{2cm}{0.5pt}\end{center}
\begin{spacing}{1}

\noindent{\sc J.~Dolbeault \& M.J.~Esteban:} CEREMADE, CNRS, UMR 7534, Universit\'e Paris-Dauphine, PSL Research University, Place de Lattre de Tassigny, F-75016 Paris, France.\\
\Email{dolbeaul@ceremade.dauphine.fr}, {\href{mailto:esteban@ceremade.dauphine.fr}{\rm\textsf{esteban@ceremade.dauphine.fr}}}

\noindent{\sc M.~Loss:} School of Mathematics, Skiles Building, Georgia Institute of Technology, Atlanta GA 30332-0160, USA.
\Email{loss@math.gatech.edu}

\end{spacing}
\end{document}